\documentclass[12pt,twoside]{article}
\usepackage{amssymb}
\font\teneufm=eufm10 scaled \magstep1
\font\seveneufm=eufm7 scaled \magstep1
\font\fiveeufm=eufm5 scaled \magstep1
\newfam\eufmfam
\textfont\eufmfam=\teneufm
\scriptfont\eufmfam=\seveneufm
\scriptscriptfont\eufmfam=\fiveeufm
\def\frak#1{{\fam\eufmfam\relax#1}}

\newfam\msbfam
\font\tenmsb=msbm10 scaled \magstep1\textfont\msbfam=\tenmsb
\font\sevenmsb=msbm7 scaled \magstep1 \scriptfont\msbfam=\sevenmsb
\font\fivemsb=msbm5 scaled \magstep1  \scriptscriptfont\msbfam=\fivemsb
\def\Bbb{\fam\msbfam \tenmsb}

\def\RR{{\Bbb R}}
\def\CC{{\Bbb C}}

\def\NN{{\Bbb N}}
\def\ZZ{{\Bbb Z}}

\def\PP{{\Bbb P}}

\def\ra{\rightarrow}

 \def\HollowBoxx #1#2#3{{\dimen0=#1 \advance\dimen0 by -#2
       \dimen1=#1 \advance\dimen1 by #3
        \vrule height 0pt depth #3 width #2
       \hskip -#3
       \vrule height #1 depth #3 width #3}}
 \def\LeftContraction{\mathord{\kern1.45pt \HollowBoxx{6pt}{3.5pt}{.4pt}}\,}

 \def\HollowBox #1#2#3{{\dimen0=#1 \advance\dimen0 by -#3
       \dimen1=#1 \advance\dimen1 by #3
        \vrule height #1 depth #3 width #3
        \vrule height 0pt depth #3 width #2
        \hskip -#3}}
 \def\RightContraction{\mathord{\, \HollowBox{6pt}{3.1pt}{.4pt}} \kern1.6pt}

\def\qed{{\hfill $\Box$}}
\newtheorem{theorem}{THEOREM}[section]

\newtheorem{proposition}[theorem]{Proposition}

\begin{document}

\begin{center}
{\Large \bf Hyperbolic
 $n$-Dimensional Manifolds
\medskip\\
with Automorphism Group 
\medskip\\
of Dimension $n^2$}
\medskip\\
\normalsize A. V. Isaev

\end{center}

\begin{quotation} \small \sl We obtain a complete classification of complex Kobayashi-hyperbolic manifolds of dimension $n\ge 2$, for which the dimension of the group of holomorphic automorphisms is equal to $n^2$. 
\end{quotation}

\thispagestyle{empty}

\pagestyle{myheadings}
\markboth{A. V. Isaev}{Automorphism Groups of Hyperbolic Manifolds}

\setcounter{section}{0}
\section{Introduction}

\setcounter{equation}{0}

For a connected complex manifold $M$ we denote by $\hbox{Aut}(M)$ the group of holomorphic automorphisms of $M$. Equipped with the natural compact-open topology, $\hbox{Aut}(M)$ is a topological group. We are interested in characterizing complex manifolds in terms of the properties of their automorphism groups.

If $M$ is Kobayashi-hyperbolic, then $\hbox{Aut}(M)$ carries the structure of a Lie group whose topology agrees with the compact-open topology \cite{Ko}, \cite{Ka}. Let $d(M):=\hbox{dim}\,\hbox{Aut}(M)$. It is known from the classical results of Kobayashi and Kaup (see \cite{Ko}, \cite{Ka}) that $d(M)\le n^2+2n$, and that $d(M)= n^2+2n$ if and only if $M$ is holomorphically equivalent to the unit ball $B^n\subset\CC^n$, where $n:=\hbox{dim}_{\CC}M$. In \cite{IKra} lower automorphism group dimensions were studied and it was shown that, for $n\ge 2$, there are in fact no hyperbolic manifolds with $n^2+3\le d(M)\le n^2+2n-1$, and that the only manifolds with $n^2<d(M)\le n^2+2$ are, up to holomorphic equivalence, $B^{n-1}\times\Delta$ (where $\Delta$ is the unit disc in $\CC$) and the 3-dimensional Siegel space (the symmetric bounded domain of type $(\hbox{III}_2)$ in $\CC^3$).

The arguments in \cite{IKra} rely on the important observation made in \cite{Ka} that $M$ is homogeneous whenever $d(M)>n^2$ and on existing techniques for homogeneous manifolds. On the other hand, if $d(M)=n^2$, the manifold $M$ need not be homogeneous. Indeed, the automorphism group of any spherical shell $S_r:=\{z\in\CC^n: r<|z|<1\}$, with $0\le r<1$, is the group $U_n$ of unitary transformations of $\CC^n$ (and thus has dimension $n^2$), but the action of $U_n$ on $S_r$ is not transitive. This simple example shows that the case $d(M)=n^2$ requires a different approach.

In \cite{GIK} we classified all hyperbolic Reinhardt domains in $\CC^n$ with $d(M)=n^2$. This classification is based on the description of the automorphism group of a hyperbolic Reinhardt domain obtained in \cite{Kru}. Further, in \cite{KV} simply connected complete hyperbolic manifolds with $d(M)=n^2$ were studied. The main result of \cite{KV} states that every such manifold is holomorphically equivalent to a Reinhardt domain and hence the classification in this case is a subset of the classification in \cite{GIK}. There are, however, examples of hyperbolic manifolds with $d(M)=n^2$ outside the class of Reinhardt domains. For instance, the factored spherical shell $S_r/\ZZ_m$, where $\ZZ_m$ is realized as a subgroup of scalar matrices in $U_n$, is not equivalent to any Reinhardt domain if $m>1$, but has the $n^2$-dimensional automorphism group $U_n/\ZZ_m$. Further, not all hyperbolic manifolds with $d(M)=n^2$ are complete and simply connected. Consider, for example, the Reinhardt domains ${\cal E}_{r,\theta}:=\Bigl\{(z',z_n)\in\CC^{n-1}\times\CC: |z'|<1,\, r(1-|z'|^2)^{\theta}<|z_n|<(1-|z'|^2)^{\theta}\Bigr\}$, with either $\theta\ge 0$, $0\le r<1$, or $\theta=-1$, $r=0$. None of them is simply connected, and complete domains only arise for either $\theta=0$, or $r=0$, $\theta>0$. At the same time, each of these domains has an $n^2$-dimensional automorphism group. This group consists of the maps
\begin{equation}
\begin{array}{lll}
z'&\mapsto&\displaystyle\frac{Az'+b}{cz'+d},\\
\vspace{0mm}&&\\
z_n&\mapsto&\displaystyle
\frac{e^{i\beta}z_n}{(cz'+d)^{2\theta}},
\end{array}\label{autgrp}
\end{equation}
where
$$
\left(\begin{array}{cc}
A& b\\
c& d
\end{array}
\right)
\in SU_{n-1,1},\quad\beta\in\RR.
$$

In this paper we obtain a complete classification of hyperbolic manifolds with $d(M)=n^2$ without any additional assumptions. We will now state our main result.
\newpage  

\begin{theorem}\label{main}\sl Let $M$ be a connected hyperbolic manifold of dimension $n\ge 2$ with $d(M)=n^2$. Then $M$ is holomorphically equivalent to one of the following manifolds:

\begin{equation}
\begin{array}{ll}
\hbox{(i)} &S_r/\ZZ_m,\,\, 0\le r<1,\,\,m\in\NN,\\
\vspace{0cm}&\\
\hbox{(ii)} & E_{\theta}:=\Bigl\{(z',z_n)\in\CC^{n-1}\times\CC: |z'|^2+|z_n|^{\theta}<1\Bigr\},\\
&\hspace{1.2cm}\theta>0,\,\theta\ne 2,\\
\vspace{0cm}&\\
\hbox{(iii)} & {\cal E}_{\theta}:=\Bigl\{(z',z_n)\in\CC^{n-1}\times\CC:|z'|<1,\, |z_n|<\\
&\hspace{1.2cm}(1-|z'|^2)^{\theta}\Bigr\},\,\theta<0,\\
\vspace{0cm}&\\ 
\hbox{(iv)} & {\cal E}_{r,\theta}=\Bigl\{(z',z_n)\in\CC^{n-1}\times\CC: |z'|<1,\, r(1-|z'|^2)^{\theta}<\\
&\hspace{1.3cm}|z_n|<(1-|z'|^2)^{\theta}\Bigr\},\,
\hbox{with either $\theta\ge 0$, $0\le r<1$},\\
&\hspace{1.3cm}\hbox{or $\theta<0$, $r=0$,}\\
\vspace{0cm}&\\
\hbox{(v)} & D_{r,\theta}:=\Bigl\{(z',z_n)\in\CC^{n-1}\times\CC: r\exp\left({\theta|z'|^2}\right)<|z_n|<\\
&\hspace{1.5cm}\exp\left({\theta|z'|^2}\right)\Bigr\},\,
\hbox{with either $\theta=1$, $0<r<1$},\\
&\hspace{1.5cm}\hbox{or $\theta=-1$, $r=0$,}\\
\vspace{0cm}&\\ 
\hbox{(vi)} & \Omega_{r,\theta}:=\Bigl\{(z',z_n)\in\CC^{n-1}\times\CC: |z'|<1,\,r(1-|z'|^2)^{\theta}<\\
&\hspace{1.5cm}\exp\left(\hbox{Re}\,z_n\right)<(1-|z'|^2)^{\theta}\Bigr\},\hbox{with either}\\
&\hspace{1.5cm}\hbox{$\theta=1$, $0\le r<1$ or $\theta=-1$, $r=0$,}\\
\vspace{0cm}&\\
 \hbox{(vii)} & {\frak S}:=\Bigl\{(z',z_n)\in\CC^{n-1}\times\CC: -1+|z'|^2<\hbox{Re}\,z_n<|z'|^2\Bigr\},\\
\vspace{0cm}&\\ 
\hbox{(viii)} & \Delta^3\,\,\, \hbox{(here $n=3$)},\\
\vspace{0cm}&\\ 
\hbox{(ix)} & B^2\times B^2\,\,\,  \hbox{(here $n=4$)}.
\end{array}\label{list}
\end{equation}

\noindent The manifolds on list (\ref{list}) are pairwise holomorphically non-equivalent.
\end{theorem}

The $n^2$-dimensional automorphism groups of manifolds (\ref{list}) are not hard to find and will explicitly appear during the course of proof of Theorem \ref{main}. Except those mentioned earlier and the well-known automorphism groups of $\Delta^3$ and $B^2\times B^2$, they are as follows: $\hbox{Aut}(E_{\theta})$ is obtained from formula (\ref{autgrp}) by replacing $\theta$ with $1/\theta$, $\hbox{Aut}(D_{0,-1})$ consists of all maps of the form 
$$
\begin{array}{lll}
z' & \mapsto & Uz'+a,\\
z_n & \mapsto &e^{i\beta}\exp\Bigl(-2\langle Uz',a\rangle-|a|^2\Bigr)z_n,
\end{array}
$$
where $U\in U_{n-1}$, $a\in\CC^{n-1}$, $\beta\in\RR$, and $\langle\cdot\,,\cdot\rangle$ denotes the inner product in $\CC^{n-1}$; $\hbox{Aut}(D_{r,1})$ is given by (\ref{gdelta}),
$\hbox{Aut}(\Omega_{r,\theta})$ is given by (\ref{autgrpcov}), $\hbox{Aut}({\frak S})$ consists of all maps of the form (\ref{thegroupsphpt}) with $\lambda=1$.

The domains $S_r$, $E_{\theta}$, ${\cal E}_{\theta}$, ${\cal E}_{r,\theta}$, $D_{r,\theta}$, $\Delta^3$ and $B^2\times B^2$ are the Reinhardt domains from the classification in \cite{GIK}. Each of the domains $\Omega_{r,\theta}$ is the universal cover of some ${\cal E}_{r',\theta'}$, namely, $\Omega_{r^{1/\theta},\,1}$ covers ${\cal E}_{r,\theta}$ for $0\le r<1$, $\theta>0$, and $\Omega_{0,-1}$ covers ${\cal E}_{0,\theta}$ for $\theta<0$. Note that the universal cover of ${\cal E}_{r,0}$ for every $0\le r<1$ is $B^{n-1}\times\Delta$ and hence has an automorphism group of dimension $n^2+2$. Observe also that ${\frak S}$ is the universal cover of the domain $D_{r,1}$ for every $0<r<1$, and that the universal cover of $D_{0,-1}$ is $B^n$ and thus has an automorphism group of dimension $n^2+2n$.

The only manifolds in (\ref{list}) that are both simply connected and complete hyperbolic are those listed under (ii), (viii) and (ix). They form the partial classification obtained in \cite{KV}. 

The paper is organized as follows. In Section \ref{proofmain1} we determine the dimensions of the orbits of the action on $M$ of $G(M):=\hbox{Aut}(M)^c$, the connected component of the identity of $\hbox{Aut}(M)$. It turns out that, unless $M$ is homogeneous, every $G(M)$-orbit is either a real or complex hypersurface in $M$. Furthermore, real hypersurface orbits are either spherical or Levi-flat, and in the latter case they are foliated by complex manifolds holomorphically equivalent to the ball $B^{n-1}$; there are at most two complex hypersurface orbits and each of them is holomorphically equivalent to $B^{n-1}$ (see Proposition \ref{main1}).

Next, in Section \ref{spherleviflat} we determine real hypersurface orbits up to $CR$-equivalence. In the spherical case there are five possible kinds of orbits (see Proposition \ref{spherorbitsprop}). In the Levi-flat case every orbit is shown to be equivalent to either $B^{n-1}\times\RR$ or $B^{n-1}\times S^1$ (see Proposition \ref{propleviflat}). The presence of an orbit of a particular kind in most cases determines $G(M)$ as a Lie group (see Propositions \ref{groupsdeterm} and \ref{propleviflat}).

In Section \ref{gluing} we prove Theorem \ref{main} in the non-homogeneous case by studying how orbits of the kinds found in Section \ref{spherleviflat} can be glued together. This gives us cases (i)-(vii) of (\ref{list}). In particular, it turns out that Levi-flat orbits equivalent to $B^{n-1}\times\RR$ cannot occur, since they lead to the domain $B^{n-1}\times\Delta$, and that there can be at most one complex hypersurface orbit.

In Section \ref{proofmain2} the homogeneous case is considered. We show that the only homogeneous manifolds that can occur are (viii) and (ix) of (\ref{list}). 

It is natural to ask what the lowest automorphism group dimension for which one can hope to obtain an explicit classification of hyperbolic manifolds is. In Section \ref{examples} we show that for $d(M)=n^2-2$ there is little hope of obtaining a reasonable classification even for the case of smoothly bounded Reinhardt domains in $\CC^2$. Namely, we construct an example of a family of pairwise holomorphically non-equivalent smoothly bounded Reinhardt domains in $\CC^2$ with automorphism group of dimension 2, parametrized by a set in $\RR^2$ that satisfies very mild conditions. In particular, there are no explicit formulas describing the domains in the family. Of course, this example is not surprising: it is expected that the automorphism groups of most Reinhardt domains in $\CC^n$ consist only of rotations in each variable and thus have dimension $n$. 

Hence beyond automorphism group dimension $n^2$, dimension $n^2-1$ is probably the only remaining candidate to investigate for the existence of a reasonable classification. We note that it follows from \cite{Kru} (see \cite{GIK}) that the automorphism group of a hyperbolic Reinhardt domain in $\CC^n$ cannot have dimension $n^2-1$.

Before proceeding, we would like to thank N. Kruzhilin and E. Lerman for useful discussions. 

\section{Dimensions of Orbits}\label{proofmain1}

\setcounter{equation}{0}

The action of $G(M)$ on $M$ is proper (see Satz 2.5 of \cite{Ka}), and therefore for every $p\in M$ its orbit $O(p):=\{f(p):f\in G(M)\}$ is a closed submanifold of $M$ and the isotropy subgroup $I_p:=\{f\in G(M): f(p)=p\}$ of $p$ is compact (see \cite{Ko}, \cite{Ka}). In this section we will obtain an initial classification of the orbits.

Let $L_p:=\{d_pf: f\in I_p\}$ be the linear isotropy subgroup, where $d_pf$ is the differential of a map $f$ at $p$. The group $L_p$ is a compact subgroup of $GL(T_p(M),\CC)$ isomorphic to $I_p$ by means of the isotropy representation
$$
\alpha_p:\, I_p\ra L_p, \quad \alpha_p(f)=d_pf
$$
(see e.g. Satz 4.3 of \cite{Ka}). We are now ready to formulate our orbit classification result. 

\begin{proposition}\label{main1}\sl Let $M$ be a connected hyperbolic manifold of dimension $n\ge 2$ with $d(M)=n^2$, and $p\in M$. Then
\smallskip\\

\noindent (i) either $M$ is homogeneous or $O(p)$ is a real or complex closed hypersurface in $M$;
\smallskip\\

\noindent (ii) if $O(p)$ is a real hypersurface, it is either Levi-flat and foliated by manifolds holomorphically equivalent to $B^{n-1}$, or spherical; there exist coordinates in $T_p(M)$ in which $L_p$ is either $U_{n-1}$ or $\ZZ_2\times U_{n-1}$, and the latter can only occur for Levi-flat orbits;
\smallskip\\

\noindent (iii) if $O(p)$ is a complex hypersurface, it is holomorphically equivalent to $B^{n-1}$, there exist coordinates in $T_p(M)$ in which $L_p=U_1\times U_{n-1}$, the subgroup $I_p':=\alpha_p^{-1}(U_1)$ is normal in $G(M)$, and the factor-group $G(M)/I_p'$ is isomorphic to $\hbox{Aut}(B^{n-1})$; there are no more than two complex hypersurface orbits in $M$.
\end{proposition}

\noindent {\bf Proof:} The proof is similar to in part that of Proposition 1.1 of \cite{IKru} (see also Satz 1.2 in \cite{Ka}). Let $V\subset T_p(M)$ be the tangent space to $O(p)$ at $p$. Clearly, $V$ is $L_p$-invariant. We assume now that $O(p)\ne M$ (and therefore $V\ne T_p(M)$) and consider the following three cases.
\smallskip\\

{\bf Case 1.} $d:=\hbox{dim}_{\CC}(V+iV)<n$.
\smallskip\\

Since $L_p$ is compact, one can introduce coordinates on $T_p(M)$ such that $L_p\subset U_n$. Further, the action of $L_p$ on $T_p(M)$
is completely reducible and the subspace $V+iV$ is invariant under this action. Hence $L_p$ can in fact be embedded in $U_{n-d}\times
U_d$. Since $\hbox{dim}\,O(p)\le 2d$, it follows that
$$
n^2\le (n-d)^2+d^2+2d,
$$
and therefore either $d=0$ or $d=n-1$. If $d=0$, then $O(p)=\{p\}$. Now Folgerung 1.10 of \cite{Ka} implies that $M$ is holomorphically equivalent to $B^n$, which is impossible, and thus $d=n-1$ (for an alternative proof of this fact see Lemma 3.1 of \cite{KV}).

We have
$$
n^2=\hbox{dim}\,L_p+\hbox{dim}\,O(p)\le n^2-2n+2+\hbox{dim}\,O(p).
$$
Hence $\hbox{dim}\,O(p)\ge 2n-2$ which implies that $\hbox{dim}\,O(p)=2d=2n-2$, and therefore $iV=V$, which means that $O(p)$ is a complex closed hypersurface in $M$. 

We have $L_p=U_1\times U_{n-1}$. The $U_{n-1}$-component of $L_p$ acts transitively on directions in $V$ and therefore, by the result of \cite{GK}, the orbit $O(p)$ is holomorphically equivalent to $B^{n-1}$. On the other hand, the $U_1$-component of $L_p$ acts trivially on $V$ and therefore the subgroup $I_p'=\alpha_p^{-1}(U_1)$ of $I_p$ corresponding to this component is the kernel of the action of $G(M)$ on $O(p)$ (this follows, for example, from Bochner's linearization theorem -- see also \cite{Ka}). Thus, $I_p'$ is normal in $G(M)$ and the factor-group $G(M)/I_p'$ acts effectively on $O(p)$. Since $G(M)/I_p'$ is connected and has dimension $n^2-1=\hbox{dim}\,\hbox{Aut}(B^{n-1})$, it is isomorphic to $\hbox{Aut}(B^{n-1})$.    
\smallskip\\

{\bf Case 2.} $T_p(M)=V+iV$ and $r:=\hbox{dim}_{\CC}(V\cap iV)>0$.
\smallskip\\

As above, $L_p$ can be embedded in $U_{n-r}\times U_r$ (we have
$r<n$). Moreover,
 $V\cap iV\ne V$ and, since $L_p$ preserves $V$, it follows that
$\hbox{dim}\,L_p<r^2+(n-r)^2$. We have $\hbox{dim}\,O(p)\le 2n-1$, and
therefore
$$
n^2<(n-r)^2+r^2+2n-1,
$$
which shows that either $r=1$, or $r=n-1$. It then follows that $\hbox{dim}\,L_p<n^2-2n+2$. Therefore, we have
$$
n^2=\hbox{dim}\,L_p+\hbox{dim}\,O(p)<n^2-2n+2+\hbox{dim}\,O(p).
$$
Hence $\hbox{dim}\,O(p)>2n-2$ and thus $\hbox{dim}\,O(p)=2n-1$. This yields that $O(p)$ is a real closed hypersurface in $M$.

Let $W$ be the orthogonal complement to $V\cap iV$ in $T_p(M)$. Clearly, $r=n-1$ and $\hbox{dim}_{\CC}\,W=1$. The group $L_p$ is a subgroup of $U_n$ and preserves $V$, $V\cap iV$, and $W$; hence it preserves the line $W\cap V$. Therefore, it can act only as $\pm\hbox{id}$ on $W$, that is, $L_p\subset\ZZ_2\times U_{n-1}$. Since $\hbox{dim}\,L_p=(n-1)^2$, we have either $L_p=U_{n-1}$, or $L_p=\ZZ_2\times U_{n-1}$. In particular, $L_p$ acts transitively on directions in $V\cap iV$. Hence $O(p)$ is either Levi-flat or strongly pseudoconvex. If $O(p)$ is strongly pseudoconvex, we have $L_p=U_{n-1}$. 

Suppose that $O(p)$ is Levi-flat. Then $O(p)$ is foliated by connected complex manifolds. Let $M_p$ be the leaf passing through $p$. Since $L_p$ acts transitively on directions in the tangent space $V\cap iV$ to $M_p$, it follows from \cite{GK} that $M_p$ is holomorphically equivalent to $B^{n-1}$. The same argument applies to any other leaf in $O(p)$. 

Assume now that $O(p)$ is strongly pseudoconvex. Since $L_p=U_{n-1}$, the dimension of the stability group of $O(p)$ at $p$ is greater than or equal to $(n-1)^2$, which implies that $p$ is an umbilic point of $O(p)$ (see e.g. \cite{EzhI}). The homogeneity of $O(p)$ now implies that $O(p)$ is spherical. 
\smallskip\\

{\bf Case 3.} $T_p(M)=V\oplus iV$.
\smallskip\\

In this case $\hbox{dim}\, V=n$ and
$L_p$ can be embedded in the real orthogonal group $O_n(\RR)$,
and therefore
$$
\hbox{dim}\,L_p+\hbox{dim}\, O(p)\le \frac{n(n-1)}{2}+n<n^2,
$$
which is a contradiction.

We will now show that there can be no more than two complex hypersurface orbits in $M$. Let $O(p)$ be a complex hypersurface for some $p\in M$. Since $L_p$ acts as $U_1$ on a complement to $V$ in $T_p(M)$ (see Case 1), there exists a neighborhood $U$ of $p$ such that for every $q\in U\setminus O(p)$ the values at $q$ of the vector fields on $M$ arising from the action of $G(M)$, span a codimension 1 subspace of $T_q(M)$. Hence there is always a real hypersurface orbit in $M$. Since the action of $G(M)$ on $M$ is proper, it follows that the orbit space $M/G(M)$ is homeomorphic to one of the following: $\RR$, $S^1$, $[0,1]$, $[0,1)$ (see \cite{M}, \cite{B-B}, \cite{AA1}, \cite{AA2}), and thus there can be no more than two complex hypersurface orbits in $M$.

The proof of the proposition is complete.\qed

\section{Real Hypersurface Orbits}\label{spherleviflat}

\setcounter{equation}{0}

In this section we will classify real hypersurface orbits up to $CR$-equivalence. We will deal with spherical orbits first.

\begin{proposition}\label{spherorbitsprop}\sl Let $M$ be a connected hyperbolic manifold of dimension $n\ge 2$ with $d(M)=n^2$. Assume that for a point $p\in M$ its orbit $O(p)$ is spherical. Then $O(p)$ is $CR$-equivalent to one of the following hypersurfaces:
\begin{equation}
\begin{array}{ll}
\hbox{(i)}& \hbox{a lens manifold ${\cal L}_m:=S^{2n-1}/\ZZ_m$ for some $m\in\NN$},\\
\hbox{(ii)} & \sigma:=\left\{(z',z_n)\in\CC^{n-1}\times\CC:\hbox{Re}\,z_n=|z'|^2\right\},\\
\hbox{(iii)} & \delta:=\hbox{$\left\{(z',z_n)\in\CC^{n-1}\times\CC: |z_n|=\exp\left(|z'|^2\right)\right\}$},\\
\hbox{(iv)} & \omega:=\left\{(z',z_n)\in\CC^{n-1}\times\CC:|z'|^2+\exp\left(\hbox{Re}\,z_n\right)=1\right\},\\
\hbox{(v)} & \varepsilon_{\alpha}:=\hbox{$\left\{(z',z_n)\in\CC^{n-1}\times\CC: |z'|^2+|z_n|^{\alpha}=1,\, z_n\ne 0\right\}$},\\
&\hspace{1.1cm}\hbox{for some $\alpha>0$}.
\end{array}\label{classificationspherorb}
\end{equation}
\end{proposition}

\noindent {\bf Proof:} For a connected Levi non-degenerate $CR$-manifold $S$ denote by $\hbox{Aut}_{CR}(S)$ the Lie group of its $CR$-automorphisms. Let $\tilde O(p)$ be the universal cover of $O(p)$. The connected component of the identity $\hbox{Aut}_{CR}(O(p))^c$ of $\hbox{Aut}_{CR}(O(p))$ acts transitively on $O(p)$ and therefore its universal cover $\widetilde{\hbox{Aut}}_{CR}(O(p))^c$ acts transitively on $\tilde O(p)$. Let $G$ be the (possibly non-closed) subgroup of $\hbox{Aut}_{CR}(\tilde O(p))$ that consists of all $CR$-automorphisms of $\tilde O(p)$ generated by this action. Observe that $G$ is a Lie group isomorphic to the factor-group of $\widetilde{\hbox{Aut}}_{CR}(O(p))^c$ by a discrete central subgroup. Let $\Gamma\subset \hbox{Aut}_{CR}(\tilde O(p))$ be the discrete subgroup whose orbits are the fibers of the covering $\tilde O(p)\ra O(p)$. The group $\Gamma$ acts freely properly discontinuously on $\tilde O(p)$, lies in the centralizer of $G$ in $\hbox{Aut}_{CR}(\tilde O(p))$ and is isomorphic to $H/H^c$, with $H=\pi^{-1}(I_p)$, where $\pi:\widetilde{\hbox{Aut}}_{CR}(O(p))^c\ra\hbox{Aut}_{CR}(O(p))^c$ is the covering map.   

The manifold $\tilde O(p)$ is spherical, and there is a local $CR$-isomorphism $\Pi$ from $\tilde O(p)$ onto a domain $D\subset S^{2n-1}$. By Proposition 1.4 of \cite{BS}, $\Pi$ is a covering map. Further, for every $f\in\hbox{Aut}_{CR}(\tilde O(p))$ there is $g\in\hbox{Aut}(D)$ such that 
\begin{equation}
g\circ \Pi=\Pi\circ f.\label{liftspher}
\end{equation}
Since $\tilde O(p)$ is homogeneous, (\ref{liftspher}) implies that $D$ is homogeneous as well, and $\hbox{dim}\,\hbox{Aut}_{CR}(\tilde O(p))=\hbox{dim}\,\hbox{Aut}_{CR}(D)$.

Clearly, $\dim\hbox{Aut}_{CR}(O(p))\ge n^2$ and therefore we have $\dim\hbox{Aut}_{CR}(D)\ge n^2$. All homogeneous domains in $S^{2n-1}$ are listed in Theorem 3.1 in \cite{BS}. It is not difficult to exclude from this list all domains with automorphism group of dimension less than $n^2$. This gives that $D$ is $CR$-equivalent to one of the following domains:
$$
\begin{array}{ll}
\hbox{(a)}& S^{2n-1},\\
\hbox{(b)}& S^{2n-1}\setminus\{\hbox{point}\},\\
\hbox{(c)}& S^{2n-1}\setminus\{z_n=0\}.
\end{array}
$$
Thus, $\tilde O(p)$ is respectively one of the following manifolds:
$$
\begin{array}{ll}
\hbox{(a)}& S^{2n-1},\\
\hbox{(b)}& \sigma,\\
\hbox{(c)}& \omega.
\end{array}
$$

If $\tilde O(p)=S^{2n-1}$, then by Proposition 5.1 of \cite{BS} the orbit $O(p)$ is $CR$-equivalent to a lens manifold as in (i) of (\ref{classificationspherorb}). 

Suppose next that $\tilde O(p)=\sigma$. The group $\hbox{Aut}_{CR}(\sigma)$ consists of all maps of the form
\begin{equation}
\begin{array}{lll}
z' & \mapsto & \lambda Uz'+a,\\
z_n & \mapsto & \lambda^2z_n+2\lambda\langle Uz',a\rangle+|a|^2+i\alpha,
\end{array}\label{thegroupsphpt}
\end{equation}
where $U\in U_{n-1}$, $a\in\CC^{n-1}$, $\lambda\in\RR^*$, $\alpha\in\RR$, and $\langle\cdot\,,\cdot\rangle$ is the inner product in $\CC^{n-1}$. It then follows that $\hbox{Aut}_{CR}(\sigma)=CU_{n-1}\ltimes N$, where $CU_{n-1}$ consists of all maps of the form (\ref{thegroupsphpt}) with $a=0$, $\alpha=0$, and $N$ is the Heisenberg group consisting of the maps of the form (\ref{thegroupsphpt}) with $U=\hbox{id}$ and $\lambda=1$.

Further, description (\ref{thegroupsphpt}) implies that $\hbox{dim}\,\hbox{Aut}_{CR}(\sigma)=n^2+1$, and therefore $n^2\le \hbox{dim}\,G\le n^2+1$. If $\hbox{dim}\,G=n^2+1$, then we have $G=\hbox{Aut}_{CR}(\sigma)^c$, and hence $\Gamma$ is a central subgroup of $\hbox{Aut}_{CR}(\sigma)^c$. Since the center of $\hbox{Aut}_{CR}(\sigma)^c$ is trivial, so is $\Gamma$. Thus, in this case $O(p)$ is $CR$-equivalent to the hypersurface $\sigma$.

Assume now that $\hbox{dim}\,G=n^2$. Since $G$ acts transitively on $\sigma$, we have $N\subset G$. Furthermore, since $G$ is of codimension 1 in $\hbox{Aut}_{CR}(\sigma)$, it contains the subgroup $SU_{n-1}\ltimes N$. By Proposition 5.6 of \cite{BS}, we have $\Gamma\subset U_{n-1}\ltimes N$. The centralizer of $SU_{n-1}\ltimes N$ in $U_{n-1}\ltimes N$ consists of all maps of the form
\begin{equation}
\begin{array}{lll}
z' &\mapsto & z',\\
z_n &\mapsto & z_n+i\alpha,
\end{array}\label{center}
\end{equation}
where $\alpha\in\RR$. Since $\Gamma$ acts freely properly discontinuously on $\sigma$, it is generated by a single map of the form (\ref{center}) with $\alpha=\alpha_0\in\RR^*$. The hypersurface $\sigma$ covers the hypersurface
\begin{equation}
\left\{(z',z_n)\in\CC^{n-1}\times\CC: |z_n|=\exp\left(\frac{2\pi}{\alpha_0}|z'|^2\right)\right\}\label{intermediate1}
\end{equation}
by means of the map
\begin{equation}
\begin{array}{lll}
z' & \mapsto & z',\\
z_n & \mapsto & \exp\left(\displaystyle\frac{2\pi}{\alpha_0} z_n\right),
\end{array}\label{coverrr}
\end{equation}
and the fibers of this map are the orbits of $\Gamma$. Hence $O(p)$ is $CR$-equivalent to hypersurface (\ref{intermediate1}). Replacing if necessary $z_n$ by $1/z_n$ we obtain that $O(p)$ is $CR$-equivalent to  the hypersurface $\delta$.                 

Suppose finally that $\tilde O(p)=\omega$. First, we will determine the group $\hbox{Aut}_{CR}(\omega)$. The general form of a $CR$-automorphism of $S^{2n-1}\setminus\{z_n=0\}$ is given by formula (\ref{autgrp}) with $\theta=1/2$ and the covering map $\Pi$ by the formula      
$$
\begin{array}{l}
z'\mapsto z',\\
z_n\mapsto \exp\left(\displaystyle\frac{z_n}{2}\right).
\end{array}
$$
Using (\ref{liftspher}) we then obtain the general form of a $CR$-automorphism of $\omega$ as follows:
\begin{equation}
\begin{array}{lll}
z'&\mapsto&\displaystyle\frac{Az'+b}{cz'+d},\\
\vspace{0mm}&&\\
z_n&\mapsto&\displaystyle z_n-2\ln(cz'+d)+i\beta,
\end{array}\label{autgrpcov}
\end{equation}
where
$$
\left(\begin{array}{cc}
A& b\\
c& d
\end{array}
\right)
\in SU_{n-1,1},\quad\beta\in\RR.
$$
In particular, $\hbox{Aut}_{CR}(\omega)$ is a connected group of dimension $n^2$ and therefore $G=\hbox{Aut}_{CR}(\omega)$. Hence $\Gamma$ is a central subgroup of $\hbox{Aut}_{CR}(\omega)$. It follows from formula (\ref{autgrpcov}) that the center of $\hbox{Aut}_{CR}(\omega)$ consists of all maps of the form (\ref{center}). Hence $\Gamma$ is generated by a single such map with $\alpha=\alpha_0\in\RR$. If $\alpha_0=0$, the orbit $O(p)$ is $CR$-equivalent to $\omega$. Let $\alpha_0\ne 0$. The hypersurface $\omega$ covers the hypersurface
\begin{equation}
\left\{(z',z_n)\in\CC^{n-1}\times\CC: |z'|^2+|z_n|^{\frac{\alpha_0}{2\pi}}=1,\, z_n\ne 0\right\}\label{intermediate2}
\end{equation}
by means of map (\ref{coverrr}). Since the fibers of this map are the orbits of $\Gamma$, it follows that $O(p)$ is $CR$-equivalent to hypersurface (\ref{intermediate2}). Replacing if necessary $z_n$ by $1/z_n$, we obtain that $O(p)$ is $CR$-equivalent to the hypersurface $\varepsilon_{\alpha}$ for some $\alpha>0$.

The proof is complete.\qed 
\smallskip\\

We will now show that in most cases the presence of a spherical orbit of a particular kind in $M$ determines the group $G(M)$ as a Lie group. Suppose that for some $p\in M$ the orbit $O(p)$ is spherical, and let ${\frak m}$ be the manifold from list (\ref{classificationspherorb}) to which $O(p)$ is $CR$-equivalent. Since $G(M)$ acts effectively on $O(p)$, the $CR$-equivalence induces an isomorphism between $G(M)$ and a (possibly non-closed) connected $n^2$-dimensional subgroup $R_{\frak m}$ of $\hbox{Aut}_{CR}({\frak m})$. A priori, $R_{\frak m}$ depends on the choice of a $CR$-isomorphism between $O(p)$ and ${\frak m}$, but, as we will see below, this dependence is insignificant.

We will now prove the following proposition.

\begin{proposition}\label{groupsdeterm} \sl${}$\linebreak

\noindent(i) $R_{S^{2n-1}}$ is conjugate to $U_n$ in $\hbox{Aut}(B^n)$,  and $R_{{\cal L}_m}=U_n/\ZZ_m$ for $m>1$; 

\noindent (ii) $R_{\sigma}$ is a subgroup in the group of all maps of the form (\ref{thegroupsphpt}) that contains $SU_{n-1}\ltimes N$;

\noindent (iii) $R_{\delta}$ consists of all maps of the form   
\begin{equation}
\begin{array}{lll}
z' & \mapsto & Uz'+a,\\
z_n & \mapsto &e^{i\beta}\exp\Bigl(2\langle Uz',a\rangle+|a|^2\Bigr)z_n,
\end{array}\label{gdelta}
\end{equation}
where $U\in U_{n-1}$, $a\in\CC^{n-1}$, $\beta\in\RR$;

\noindent (iv) $R_{\omega}$ consists of all maps of the form (\ref{autgrpcov});

\noindent (v) $R_{\varepsilon_{\alpha}}$ consists of all maps of the form (\ref{autgrp}) with $\theta=1/\alpha$.
\end{proposition}

\noindent{\bf Proof:} Suppose first that ${\frak m}={\cal L}_m$, for some $m\in\NN$. Then $O(p)$ is compact and, since $I_p$ is compact as well, it follows that $G(M)$ is compact. Assume first that $m=1$. In this case $R_{S^{2n-1}}$ is a subgroup of $\hbox{Aut}_{CR}(S^{2n-1})=\hbox{Aut}(B^n)$.
Since $R_{S^{2n-1}}$ is compact, it is conjugate to a subgroup of $U_n$, which is a maximal compact subgroup in $\hbox{Aut}(B^n)$. Since both $R_{S^{2n-1}}$ and $U_n$ are $n^2$-dimensional, $R_{S^{2n-1}}$ is conjugate to the full group $U_n$. Suppose now that $m>1$. It is straightforward to determine the group $\hbox{Aut}_{CR}\left({\cal L}_m\right)$ by lifting $CR$-automorphisms of ${\cal L}_m$ to its universal cover $S^{2n-1}$. This group is $U_n/\ZZ_m$ acting on $\CC^n\setminus\{0\}/\ZZ_m$ in the standard way. In particular, $\hbox{Aut}_{CR}\left({\cal L}_m\right)$ is connected and has dimension $n^2$. Therefore, $R_{{\cal L}_m}=U_n/\ZZ_m$.

Assume now that ${\frak m}=\sigma$. The group $\hbox{Aut}_{CR}(\sigma)$ consists of all maps of the form (\ref{thegroupsphpt}) and has dimension $n^2+1$. 
Since $R_{\sigma}$ acts transitively on $\sigma$, it contains the subgroup $N$ (see the proof of Proposition \ref{spherorbitsprop}). Furthermore, $R_{\sigma}$ is a codimension 1 subgroup of $\hbox{Aut}_{CR}(\sigma)$, and thus contains the subgroup $SU_{n-1}\ltimes N$.

Next, the groups $\hbox{Aut}_{CR}(\delta)$, $\hbox{Aut}_{CR}(\omega)$, $\hbox{Aut}_{CR}(\varepsilon_{\alpha})$ are $n^2$-dimensional and connected. Indeed, $\hbox{Aut}_{CR}(\delta)$ can be determined by considering the universal cover of $\delta$ (see the proof of Proposition \ref{spherorbitsprop}) and consists of all maps of the form (\ref{gdelta}). The group $\hbox{Aut}_{CR}(\omega)$ was found in the proof of Proposition \ref{spherorbitsprop} and consists of all maps of the form (\ref{autgrpcov}). The group $\hbox{Aut}_{CR}(\varepsilon_{\alpha})$ can be found by considering the universal cover of $\varepsilon_{\alpha}$ and consists of all maps of the form (\ref{autgrp}) with $\theta=1/\alpha$. Thus, we obtain statements (iii), (iv) and (v) of the proposition. 

The proof is complete.\qed
\smallskip\\

Next, we will classify Levi-flat orbits.

\begin{proposition}\label{propleviflat} \sl Let $M$ be a connected complex hyperbolic manifold of dimension $n\ge 2$ with $d(M)=n^2$. Assume that for a point $p\in M$ its orbit $O(p)$ is Levi-flat. Then $O(p)$ is equivalent to either $B^{n-1}\times\RR$ or $B^{n-1}\times S^1$ by means of a real-analytic $CR$-map. The $CR$-equivalence can be chosen so that it transforms $G(M)$ into the group $R_{B^{n-1}\times\RR}$ of all maps of the form
$$
(z',z_n)\mapsto (a(z'),z_n+b),
$$
in the first case, and into the group $R_{B^{n-1}\times S^1}$ of all maps of the form
$$
(z',z_n)\mapsto (a(z'),e^{ic}z_n),
$$
in the second case, where $a\in\hbox{Aut}(B^{n-1})$, $b,c\in\RR$.
\end{proposition}

\noindent{\bf Proof:} Recall that the hypersurface $O(p)$ is foliated by complex manifolds holomorphically equivalent to $B^{n-1}$ (see (ii) of Proposition \ref{main1}). Denote by ${\frak g}$ the Lie algebra of vector fields on $O(p)$ arising from the action of $G(M)$. Clearly, ${\frak g}$ is isomorphic to the Lie algebra of $G(M)$. For $q\in O(p)$ we consider the leaf $M_q$ of the foliation passing through $q$ and the subspace ${\frak l}_q\subset{\frak g}$ of all vector fields tangent to $M_q$ at $q$. Since vector fields in ${\frak l}_q$ remain tangent to $M_q$ at each point in $M_q$, the subspace ${\frak l}_q$ is in fact a Lie subalgebra of ${\frak g}$. It follows from the definition of ${\frak l}_q$ that $\hbox{dim}\,{\frak l}_q=n^2-1$.

Denote by $H_q$ the (possibly non-closed) connected subgroup of $G(M)$ with Lie algebra ${\frak l}_q$. It is straightforward to verify that the group $H_q$ acts on $M_q$ by holomorphic transformations and that $I_q^c\subset H_q$. If some element $g\in H_q$ acts trivially on $M_q$, then $g\in I_q$. If $I_q$ is isomorphic to $U_{n-1}$, every element of $I_q$ acts non-trivially on $M_q$ and thus $g=\hbox{id}$; if $I_q$ is isomorphic to $\ZZ_2\times U_{n-1}$ and $g\ne\hbox{id}$, then $g=g_q$, where $g_q$ denotes the element of $I_q$ corresponding to the non-trivial element in $\ZZ_2$ (see Case 2 in the proof of Proposition \ref{main1}). Thus, either $H_q$ or $H_q/\ZZ_2$ acts effectively on $M_q$ (the former case occurs if $g_q\not\in H_q$, the latter if $g_q\in H_q$). Since $\hbox{dim}\,H_q=n^2-1=\hbox{dim}\,\hbox{Aut}(B^{n-1})$, we obtain that either $H_q$ or $H_q/\ZZ_2$ is isomorphic to $\hbox{Aut}(B^{n-1})$.

We will now show that in fact $g_q\not\in H_q$. Assuming the opposite, we have $I_q\subset H_q$. It then follows that $I_q$ is a maximal compact subgroup of $H_q$ since its image under the projection $H_q\ra \hbox{Aut}(B^{n-1})$ is a maximal compact subgroup of $\hbox{Aut}(B^{n-1})$. However, every maximal compact subgroup of a connected Lie group is connected whereas $I_q$ is not. Thus, $g_q\not\in H_q$, and hence $H_q$ is isomorphic to $\hbox{Aut}(B^{n-1})$.     

We will now show that $H_{q_1}=H_{q_2}$ for all $q_1,q_2\in O(p)$. Suppose first that $n\ge 3$. The Lie algebra ${\frak l}_q$ is isomorphic to ${\frak{su}}_{n-1,1}$. Since the algebra ${\frak{su}}_{n-1,1}$ does not have codimension 1 subalgebras if $n\ge 3$ (see e.g. \cite{EaI}), ${\frak l}_q$ is the only codimension 1 subalgebra of ${\frak g}$ in this case. This implies that $H_{q_1}=H_{q_2}$ for all $q_1,q_2\in O(p)$ for $n\ge 3$. Assume now that $n=2$. By Bochner's theorem there exist a local holomorphic change of coordinates $F$ near $q$ on $M$ that identifies an $I_q^c$-invariant neighborhood $U$ of $q$ with an $L_q^c$-invariant neighborhood $U'$ of the origin in $T_q(M)$ such that $F(q)=0$ and $F(gs)=\alpha_q(g)F(s)$ for all $g\in I_q^c$ and $s\in U$ (here $L_q$ is the linear isotropy group and $\alpha_q$ is the isotropy representation at $q$). In the proof of Proposition \ref{main1} (see Case 2) we have seen that $L_q^c$ fixes every point in the orthogonal complement $W_q$ to $T_q(M_q)$ in $T_q(M)$. Therefore, for all points $s$ in the curve $F^{-1}\Bigl(U'\cap\Bigl(W_q\cap T_q(O(p))\Bigr)\Bigr)$ we have $I_s^c=I_q^c$. In particular, $H_q\cap H_s$ contains a subgroup isomorphic to $U_1$ for such $s$. Hence ${\frak l}_q\cap{\frak l}_s$ is isomorphic to a subalgebra of ${\frak {su}}_{1,1}$ of dimension at least 2 that contains the subalgebra of all diagonal matrices in ${\frak {su}}_{1,1}$. It is straightforward to check, however, that every such subalgebra must coincide with all of ${\frak {su}}_{1,1}$, which shows that $H_q=H_s$ for all $s\in F^{-1}\Bigl(U'\cap\Bigl(W_q\cap T_q(O(p))\Bigr)\Bigr)$ and hence for all $s$ in a neighborhood of $q$. Since this argument can be applied to any $q\in O(p)$, we obtain that $H_{q_1}=H_{q_2}$ for all $q_1,q_2\in O(p)$ if $n=2$ as well. From now on we denote the coinciding groups $H_q$ by $H$ and the coinciding algebras ${\frak l}_q$ by ${\frak l}$. Note that the above argument also shows that $H$ is a normal subgroup of $G(M)$ and ${\frak l}$ is an ideal in ${\frak g}$.

We will now show that $H$ is closed in $G(M)$. Let $U$ be a neighborhood of $0$ in ${\frak g}$ where the exponential map into $G(M)$ is a diffeomorphism, and let $V:=\exp(U)$. To prove that $H$ is closed in $G(M)$ it is sufficient to show that for some neighborhood $W$ of $\hbox{id}\in G(M)$, $W\subset V$, we have $H\cap W= \exp({\frak l}\cap U)\cap W$. Assuming the opposite we obtain a sequence $\{h_j\}$ of elements of $H$ converging to $\hbox{id}$ in $G(M)$ such that for every $j$ we have $h_j=\exp(a_j)$ with $a_j\in U\setminus{\frak l}$. Fix $q\in O(p)$. There exists a neighborhood ${\cal V}$ of $q$ in $O(p)$ that is $CR$-equivalent to the direct product of a segment in $\RR$ and a ball in $\CC^{n-1}$. For every $s$ in this neighborhood we denote by $N_s$ the complex hypersurface lying in ${\cal V}$ that arises from this representation and passes through $s$. We call such hypersurfaces local leaves. Let $q_j:=h_jq$. If $j$ is sufficiently large, $q_j\in{\cal V}$. For every $j$, the local leaf $N_{q_j}$ is distinct from $N_q$, and $N_{q_j}$ accumulate to $N_q$. At the same time we have $N_q\subset M_q$ and $N_{q_j}\subset M_q$ for all $j$ and thus the leaf $M_q$ accumulates to itself. Below we will show that this is in fact impossible thus obtaining a contradiction.

Consider the action of $I_q^c$ on $M_q$. The orbit $O'(s)$ of every point $s\in M_q$, $s\ne q$, of this action is diffeomorphic to the sphere $S^{2n-3}$ and we have $M_q\setminus O'(s)=V_1(s)\cup V_2(s)$, where $V_1(s)$, $V_2(s)$ are open in $M_q$, disjoint, the point $q$ lies in $V_1(s)$, and $V_2(s)$ is diffeomorphic to a spherical shell in $\CC^{n-1}$. Since $I_q^c$ is compact, there exist neighborhoods ${\cal V}$, ${\cal V}'$ of $q$ in $O(p)$, ${\cal V}'\subset{\cal V}$, such that ${\cal V}$ is represented as a union of local leaves, and $I_q^c{\cal V}'\subset {\cal V}$. Since $M_q$ accumulates to itself near $q$, there exists $s_0\in M_q\cap{\cal V}'$, $s_0\not\in N_q$. Clearly, $O'(s_0)\subset{\cal V}$. Since the vector fields in ${\frak g}$ arising from the action of $I_q^c$ on $O(p)$ are tangent to $N_{s_0}\subset M_q$ at $s_0$ and ${\cal V}$ is partitioned into non-intersecting local leaves, the orbit $O'(s_0)$ lies in $N_{s_0}$. Then we have $N_{s_0}\setminus O'(s_0)=W_1(s_0)\cup W_2(s_0)$, where $W_1(s_0)$ is diffeomorphic to the ball. Since $q\not\in N_{s_0}$, we have and $V_2(s_0)=W_1(s_0)$, which is impossible since in this case $V_2(s_0)$ is not diffeomorphic to a spherical shell. This contradiction shows that $H$ is closed in $G(M)$.

Thus, since the action of $G(M)$ is proper on $M$, the action of $H$ on $O(p)$ is proper as well, and the orbits of this action are the leaves of the foliation. Since all isotropy subgroups for the action of $H$ on $O(p)$ are conjugate to each other, $O(p)/H$ is homeomorphic to either $\RR$ or $S^1$ (see \cite{M}, \cite{B-B}, \cite{AA1}, \cite{AA2}).
By \cite{A}, \cite{P} there is a complete $H$-invariant Riemannian metric on $O(p)$. Fix $p_0\in O(p)$, let $F:M_{p_0}\ra B^{n-1}$ be a biholomorphism, and consider a normal geodesic $\gamma_0$ emanating from $p_0$.

Suppose first that $O(p)/H$ is homeomorphic to $\RR$. Then $\gamma_0$ is diffeomorphic to $\RR$ (see \cite{AA1}). We will now construct a $CR$-isomorphism $F_1:O(p)\ra B^{n-1}\times\RR$ using the properties of normal geodesics listed in Proposition 4.1 of \cite{AA1}. For $q\in O(p)$ consider $M_q$ and let $r$ be the (unique) point where $\gamma_0$ intersects $M_q$. Let $h\in H$ be such that $q=hr$ and $t_q$ be the value of the parameter on $\gamma$ such that $\gamma(t_q)=r$. Then we set $F_1(q):=(F(hp_0),t_q)$. By construction, $F_1$ is a real-analytic $CR$-map.  

Suppose now that $O(p)/H$ is homeomorphic to $S^1$. Following \cite{AA1}, consider the Weyl group $W(\gamma_0)$ of $\gamma_0$.
This group can be identified with a subgroup of the group $N_H(I_q^c)/I_q^c$, where $q$ is any point in $\gamma_0$ and $N_H(I_q^c)$ is the normalizer of $I_q^c$ in $H$. Since $H$ is isomorphic to $\hbox{Aut}(B^{n-1})$ and $I_q^c$ upon this identification is conjugate to $U_{n-1}\subset\hbox{Aut}(B^{n-1})$, we see that $N_H(I_q^c)=I_q^c$ and thus $W(\gamma_0)$ is trivial. This implies that $\gamma_0$ is diffeomorphic to $S^1$ (see \cite{AA1}). We will now construct a $CR$-isomorphism $F_2:O(p)\ra B^{n-1}\times S^1$, using, as before, the properties of normal geodesics from Proposition 4.1 of \cite{AA1}. For $q\in O(p)$ consider $M_q$ and let $r$ be the point where $\gamma_0$ intersects $M_q$. Let $h\in H$ be such that $q=hr$ and $t_q$ be the least value of the parameter on $\gamma_0$ such that $\gamma_0(t_q)=r$. Then we set $F_2(q):=(F(hp_0),e^{2\pi it_q/T})$, where $T$ is the least positive value of the parameter on $\gamma_0$ such that $\gamma(T)=p_0$. The $CR$-map $F_2$ is real-analytic.

It follows from the construction of the maps $F_j$ that the elements of $H':=F_j\circ H|_{O(p)}\circ F_j^{-1}$ have the form
$$
(z',u)\mapsto (a(z'),u),
$$
for $j=1,2$, where $a\in\hbox{Aut}(B^{n-1})$ and $(z',u)$ is a point in either $B^{n-1}\times\RR$, or $B^{n-1}\times S^1$, respectively. We will now find the general form of the elements of the group $G':=F_j\circ G(M)|_{O(p)}\circ F_j^{-1}$. 

Every $CR$-isomorphism of either $B^{n-1}\times\RR$, or $B^{n-1}\times S^1$ has the form
\begin{equation}
(z',u)\mapsto (a_u(z'),\mu(u)),\label{prodgeneral}
\end{equation}
where $a_u\in\hbox{Aut}(B^{n-1})$ for every $u$ and $\mu$ is a diffeomorphism of either $\RR$ or $S^1$, respectively. Since $H$ is normal in $G(M)$, $H'$ is normal in $G'$. Fix an element of $G'$ and let $\{a_u\}$ be the corresponding family of automorphisms of $B^{n-1}$, as in (\ref{prodgeneral}). Then we have $a_{u_1} a a_{u_1}^{-1}=a_{u_2} a a_{u_2}^{-1}$ for all $a\in\hbox{Aut}(B^{n-1})$ and all $u_1, u_2$. Therefore, $a_{u_1}a_{u_2}^{-1}$ lies in the center of $\hbox{Aut}(B^{n-1})$, which is trivial. Hence we obtain that $a_{u_1}=a_{u_2}$ for all $u_1,u_2$. This shows that every element of $G'$ is a composition of an element of $H'$ and an element of a one-parameter family of real-analytic automorphism of the form
$$
(z',u)\mapsto (z',\mu(u)).
$$

Since for every $q\in O(p)$ there can exist at most one element of $G(M)$ outside $H$ that fixes $q$ (namely, the transformation $g_q$), the corresponding one-parameter family $\{\mu_{\tau}\}_{\tau\in\RR}$ has no fixed points in either $\RR$ or $S^1$, respectively. If $O(p)$ is equivalent to $B^{n-1}\times\RR$, then under the diffeomorphism of $\RR$ inverse to the map $x\mapsto\mu_x(0)$, the family $\{\mu_{\tau}\}_{\tau\in\RR}$ transforms into the family  
$$
(z',u)\mapsto (z',u+\tau).
$$
If $O(p)$ is equivalent to $B^{n-1}\times S^1$, then the diffeomorphism of $S^1$
$$
e^{i\theta}\mapsto\exp\left(\frac{2\pi i}{c}\int_{0}^{\theta}\frac{1}{V(t)}\,dt\right),
$$
with $V(t):=d\mu_{\tau}/d\tau|_{\tau=0,u=e^{it}}$ and $c=\int_{0}^{2\pi}dt/V(t)$, $0\le t,\theta\le 2\pi$, transforms the family $\{\mu_{\tau}\}_{\tau\in\RR}$ into the family 
$$
(z',u)\mapsto (z',e^{\frac{2\pi i}{c}\tau}u).
$$
The diffeomorphisms of $\RR$ and $S^1$ constructed above are real-analytic. 

The proof is complete. \qed
\smallskip\\

The hypersurfaces determined in Propositions \ref{spherorbitsprop} and \ref{propleviflat} will be called the {\it models}\, for real hypersurface orbits. For Levi-flat orbits we will always assume that the $CR$-equivalence between the orbit and the model is chosen to ensure that $G(M)$ transforms into either $R_{B^{n-1}\times\RR}$ or $R_{B^{n-1}\times S^1}$.  

It is straightforward to see that all models are pairwise $CR$ non-equivalent. For all models, except the family $\{\varepsilon_{\alpha}\}_{\alpha>0}$, this follows from topological considerations and comparison of the automorphism groups. Furthermore, any $CR$-isomorphism between $\varepsilon_{\alpha}$ and $\varepsilon_{\beta}$ can be lifted to a $CR$-automorphism of $\omega$, the universal cover of each of $\varepsilon_{\alpha}$ and $\varepsilon_{\beta}$, and it is straightforward to observe that in order for the lifted map to be well-defined and one-to-one, we must have $\beta=\alpha$.

\section{The Non-Homogeneous Case}\label{gluing}

\setcounter{equation}{0}

In this section we will prove Theorem \ref{main} in the non-homogeneous case by studying how real and complex hypersurface orbits can be glued together to form hyperbolic manifolds with $n^2$-dimensional automorphism groups.

Suppose that for some $p\in M$ the orbit $O(p)$ is $CR$-equivalent to a lens manifold ${\cal L}_m$, for some $m\in\NN$. In this case $G(M)$ is compact, hence there are no complex hypersurface orbits and the model for every orbit is a lens manifold. Assume first that $m=1$. Then $M$ admits an effective action of $U_n$ by holomorphic transformations and therefore is holomorphically equivalent to one of the manifolds listed in \cite{IKru}. The only hyperbolic manifolds on the list with $n^2$-dimensional automorphism group are factored spherical shells $S_r/\ZZ_l$, for $0\le r<1$ and $l=|nk+1|$ with $k\in\ZZ$ (note that for such $l$ the groups $U_n$ and $U_n/\ZZ_l$ are isomorphic). However, for $l>1$ no orbit of the action of $\hbox{Aut}(S_r/\ZZ_l)=U_n/\ZZ_l$ on $S_r/\ZZ_l$ is diffeomorphic to $S^{2n-1}$, which gives that $M$ is in fact equivalent to $S_r$.

Assume now that $m>1$. Let $f: O(p)\ra {\cal L}_m$ be a $CR$-isomorphism. Then we have 
\begin{equation}
f(gq)=\varphi(g)f(q),\label{equivar}
\end{equation}
where $q\in O(p)$, for some Lie group isomorphism $\varphi: G(M)\ra U_n/\ZZ_m$. The $CR$-isomorphism $f$ extends to a biholomorphic map from a neighborhood $U$ of $O(p)$ in $M$ onto a neighborhood $W$ of ${\cal L}_m$ in $\CC^n\setminus\{0\}/\ZZ_m$. Since $G(M)$ is compact, one can choose $U$ to be a connected union of $G(M)$-orbits. Then property (\ref{equivar}) holds for the extended map, and therefore every $G(M)$-orbit in $U$ is taken onto a $U_n/\ZZ_m$-orbit in $\CC^n\setminus\{0\}/\ZZ_m$ by this map. Thus, $W=S_r^R/\ZZ_m$ for some $0\le r<R<\infty$, where $S_r^R:=\left\{z\in\CC^n: r<|z|<R\right\}$ is a spherical shell.

Let $D$ be a maximal domain in $M$ such that there exists a biholomorphic map $f$ from $D$ onto $S_r^R/\ZZ_m$ for some $r,R$, satisfying (\ref{equivar}) for all $g\in G(M)$ and $q\in D$. As was shown above, such a domain $D$ exists. Assume that $D\ne M$ and let $x$ be a boundary point of $D$. Consider the orbit $O(x)$. Let ${\cal L}_k$ for some $k>1$ be the model for $O(x)$ and $f_1:O(x)\ra{\cal L}_k$ a $CR$-isomorphism satisfying (\ref{equivar}) for $g\in G(M)$, $q\in O(x)$ and an isomorphism $\varphi_1:G(M)\ra U_n/\ZZ_k$ in place of $\varphi$. The map $f_1$ can be holomorphically extended to a neighborhood $V$ of $O(x)$ that one can choose to be a connected union of $G(M)$-orbits. The extended map satisfies (\ref{equivar}) for $g\in G(M)$, $q\in V$ and $\varphi_1$ in place of $\varphi$. For $s\in V\cap D$ we consider the orbit $O(s)$. The maps $f$ and $f_1$ take $O(s)$ into some surfaces $r_1S^{2n-1}/\ZZ_m$ and $r_2S^{2n-1}/\ZZ_k$, respectively, with $r_1,r_2>0$.
Hence $F:=f_1\circ f^{-1}$ maps $r_1S^{2n-1}/\ZZ_m$ onto $r_2S^{2n-1}/\ZZ_k$. Since ${\cal L}_m$ and ${\cal L}_k$ are not $CR$-equivalent for distinct $m$, $k$, we obtain $k=m$. Furthermore, every $CR$-isomorphism between $r_1S^{2n-1}/\ZZ_m$ and $r_2S^{2n-1}/\ZZ_m$ has the form $[z]\mapsto [r_2/r_1Uz]$, where $U\in U_n$, and $[z]\in\CC^n\setminus\{0\}/\ZZ_m$ denotes the equivalence class of a point $z\in\CC^n\setminus\{0\}$. Therefore, $F$ extends to a holomorphic automorphism of $\CC^n\setminus\{0\}/\ZZ_m$.

We claim that $V$ can be chosen so that $D\cap V$ is connected and\linebreak $V\setminus(D\cup O(x))\ne\emptyset$. Indeed, since $O(x)$ is strongly pseudoconvex and closed in $M$, for $V$ small enough we have $V=V_1\cup V_2\cup O(x)$, where $V_j$ are open connected non-intersecting sets. For each $j$, $D\cap V_j$ is a union of $G(M)$-orbits and therefore is mapped by $f$ onto a union of the quotients of some spherical shells. If there are more than one such factored shells, then there is a factored shell such that the closure of its inverse image under $f$ is disjoint from $O(x)$, and hence $D$ is disconnected which contradicts the definition of $D$. Thus, $D\cap V_j$ is connected for $j=1,2$, and, if $V$ is sufficiently small, then each $V_j$ is either a subset of $D$ or disjoint from it. If $V_j\subset D$ for $j=1,2$, then $M=D\cup V$ is compact, which is impossible since $M$ is hyperbolic and $d(M)>0$. Therefore, for some $V$ there is only one $j$ for which $D\cap V_j\ne\emptyset$. Thus, $D\cap V$ is connected and $V\setminus(D\cup O(x))\ne\emptyset$, as required. 

Setting now
\begin{equation}
\tilde f:=\Biggl\{\begin{array}{l}
f\hspace{1.7cm}\hbox{on $D$}\\
F^{-1}\circ f_1\hspace{0.45cm}\hbox{on $V$},
\end{array}\label{extens}
\end{equation}
we obtain a biholomorphic extension of $f$ to $D\cup V$. By construction, $\tilde f$ satisfies (\ref{equivar}) for $g\in G(M)$ and $q\in D\cup V$. Since $D\cup V$ is strictly larger than $D$, we obtain a contradiction with the maximality of $D$. Thus, we have shown that in fact $D=M$, and hence $M$ is holomorphically equivalent to $S_{r/R}/\ZZ_m$.

We will now consider the case when all orbits in $M$ are non-compact. First, we will assume that every orbit in $M$ is a real hypersurface. In this situation we will use an orbit gluing procedure illustrated above for the case of compact orbits. It comprises the following steps:
\vspace{0cm}\\

\noindent (1). Start with a real hypersurface orbit $O(p)$ with model ${\frak m}$ and consider a real-analytic $CR$-isomorphism $f:O(p)\ra{\frak m}$ that satisfies (\ref{equivar}) for all $g\in G(M)$ and $q\in O(p)$, where $\varphi:G(M)\ra R_{\frak m}$ is a Lie group isomorphism. 
\vspace{0cm}\\

\noindent (2). Verify that for every model ${\frak m}'$ the group $R_{{\frak m}'}$ acts by holomorphic transformations with real hypersurface orbits on a domain ${\cal D}\subset\CC^n$ containing ${\frak m}'$. Observe that for every model ${\frak m}'$ every orbit of the action of $R_{{\frak m}'}$ on ${\cal D}$ is a real hypersurface $CR$-equivalent to ${\frak m}'$.
\vspace{0cm}\\ 

\noindent (3). Observe that $f$ can be extended to a biholomorphic map from a $G(M)$-invariant connected neighborhood of $O(p)$ in $M$ onto an $R_{\frak m}$-invariant neighborhood of ${\frak m}$ in ${\cal D}$. First of all, extend $f$ to some neighborhood $U$ of $O(p)$ to a biholomorphic map onto a neighborhood $W$ of ${\frak m}$ in $\CC^n$. Let $W'=W\cap {\cal D}$ and $U'=f^{-1}(W')$. Fix $s\in U'$ and $s_0\in O(s)$. Choose $h_0\in G(M)$ such that $s_0=h_0s$ and define $f(s_0):=\varphi(h_0)f(s)$. To see that $f$ is well-defined at $s_0$, suppose that for some $h_1\in G(M)$, $h_1\ne h_0$, we have $s_0=h_1s$, and show that $\varphi(h)$ fixes $f(s)$, where $h:=h_1^{-1}h_0$. Indeed, for every $g\in G(M)$ identity (\ref{equivar}) holds for $q\in U_g$, where $U_g$ is the connected component of $g^{-1}(U')\cap U'$ containing $O(p)$. Since $h\in I_s$, we have $s\in U_h$ and the application of (\ref{equivar}) to $h$ and $s$ yields that $\varphi(h)$ fixes $f(s)$, as required. Thus, $f$ extends to $U'':=\cup_{q\in U'}O(q)$. The extended map satisfies (\ref{equivar}) for all $g\in G(M)$ and $q\in U''$.
\vspace{0cm}\\

\noindent (4). Consider a maximal $G(M)$-invariant domain $D\subset M$ from which there exists a biholomorphic map $f$ onto an $R_{\frak m}$-invariant domain in ${\cal D}$ satisfying (\ref{equivar}) for all $g\in G(M)$ and $q\in D$. The existence of such a domain is guaranteed by the previous step. Assume that $D\ne M$ and consider $x\in\partial D$. Let ${\frak m}_1$ be the model for $O(x)$ and let $f_1:O(x)\ra{\frak m}_1$ be a real-analytic $CR$-isomorphism satisfying (\ref{equivar}) for all $g\in G(M)$, $q\in O(x)$ and some Lie group isomorphism $\varphi_1: G(M)\ra R_{{\frak m}_1}$ in place of $\varphi$. 
Let ${\cal D}_1$ be the domain in $\CC^n$ containing ${\frak m}_1$ on which $R_{{\frak m}_1}$ acts by holomorphic transformations with real hypersurface orbits $CR$-equivalent to ${\frak m}_1$. As in (3), extend $f_1$ to a biholomorphic map from a connected $G(M)$-invariant neighborhood $V$ of $O(x)$ onto an $R_{{\frak m}_1}$-invariant neighborhood of ${\frak m}_1$ in ${\cal D}_1$. The extended map satisfies (\ref{equivar}) for all $g\in G(M)$, $q\in V$ and $\varphi_1$ in place of $\varphi$. Consider $s\in V\cap D$. The maps $f$ and $f_1$ take $O(s)$ onto an $R_{\frak m}$-orbit in ${\cal D}$ and an $R_{{\frak m}_1}$-orbit in ${\cal D}_1$, respectively. Then $F:=f_1\circ f^{-1}$ maps the $R_{\frak m}$-orbit onto the $R_{{\frak m}_1}$-orbit. Since all models are pairwise $CR$ non-equivalent, we obtain ${\frak m}_1={\frak m}$. 
\vspace{0cm}\\   

\noindent (5). Show that $F$ extends to a holomorphic automorphism of ${\cal D}$. For spherical ${\frak m}$ this will follow from the fact that $F$ maps an $R_{\frak m}$-orbit onto an $R_{\frak m}$-orbit, for Levi-flat ${\frak m}$ a slightly more detailed analysis will be required.
\vspace{0cm}\\ 

\noindent (6). Show that $V$ can be chosen so that $D\cap V$ is connected and\linebreak $V\setminus(D\cup O(x))\ne\emptyset$. This follows from the hyperbolicity of $M$ and the existence of a neighborhood $V'$ of $O(x)$ such that $V'=V_1\cup V_2\cup O(x)$, where $V_j$ are open connected non-intersecting sets. For spherical ${\frak m}$ the existence of such $V'$ follows, for example, from the strong pseudoconvexity of ${\frak m}$, for Levi-flat ${\frak m}$ it follows from the explicit form of the models: indeed, each of $B^{n-1}\times\RR$, $B^{n-1}\times S^1$ splits $B^{n-1}\times\CC$.
\vspace{0cm}\\ 

\noindent (7). Use formula (\ref{extens}) to extend $f$ to $D\cup V$ thus obtaining a contradiction with the maximality of $D$. This shows that in fact $D=M$ and hence $M$ is biholomorphically equivalent to an $R_{\frak m}$-invariant domain in ${\cal D}$. In all the cases below the determination of $R_{\frak m}$-invariant domains will be straightforward, and a classification of manifolds $M$ not containing complex hypersurface orbits will follow.
\vspace{0cm}\\ 

We will now apply our general gluing procedure to each model using the notation introduced above as well as in Propositions \ref{spherorbitsprop}, \ref{groupsdeterm} and \ref{propleviflat}. Suppose first that ${\frak m}=\sigma$. Clearly, $R_{\sigma}$ acts by holomorphic transformations on all of $\CC^n$, so in this case ${\cal D}=\CC^n$. Let $f:O(p)\ra{\frak m}$ be a real-analytic $CR$-isomorphism that satisfies (\ref{equivar}) for all $g\in G(M)$ and $q\in O(p)$, where $\varphi:G(M)\ra R_{\frak \sigma}$ is a Lie group isomorphism. As we showed in (3), the map $f$ can be extended to a biholomorphic map between a $G(M)$-invariant neighborhood $U$ of $O(p)$ in $M$ and a $R_{\sigma}$-invariant neighborhood $W$ of $\sigma$ in $\CC^n$ satisfying (\ref{equivar}) for all $g\in G(M)$ and $q\in U$. This implies that the $R_{\sigma}$-orbit of every point in $W$ is a real hypersurface in $\CC^n$, which can only happen if $\lambda=1$ for every element of $R_{\sigma}$ (see formula (\ref{thegroupsphpt})). Therefore, $R_{\sigma}=U_{n-1}\ltimes N$, and thus $R_{\sigma}$ acts with real hypersurface orbits on all of $\CC^n$. The $R_{\sigma}$-orbit of every point in $\CC^n$ is of the form
$$
\left\{(z',z_n)\in\CC^{n-1}\times\CC:\hbox{Re}\,z_n=|z'|^2+r\right\},
$$
where $r\in\RR$, and every $R_{\sigma}$-invariant domain in $\CC^n$ is given by
$$
{\frak S}_r^R:=\left\{(z',z_n)\in\CC^{n-1}\times\CC: r+|z'|^2<\hbox{Re}\,z_n<R+|z'|^2\right\},
$$
where $-\infty\le r<R\le\infty$. Every $CR$-isomorphism between two $R_{\sigma}$-orbits is a composition of a map of the form (\ref{thegroupsphpt}) and a translation in the $z_n$-variable. Therefore, $F$ in this case extends to a holomorphic automorphism of $\CC^n$. Now our gluing procedure implies that $M$ is holomorphically equivalent to ${\frak S}_r^R$ for some $-\infty\le r<R\le\infty$. Therefore, $M$ is holomorphically equivalent either to the domain ${\frak S}$ or (for $R=\infty$) to $B^n$; the latter is clearly impossible.

Assume next that ${\frak m}=\delta$. Again, we have ${\cal D}=\CC^n$. The $R_{\delta}$-orbit of every point in $\CC^n$ has the form 
$$
\left\{(z',z_n)\in\CC^{n-1}\times\CC: |z_n|=r\exp\left(|z'|^2\right)\right\},
$$
where $r>0$, and hence every $R_{\delta}$-invariant domain in $\CC^n$ is given by
$$
D_r^R:=\Bigl\{(z',z_n)\in\CC^{n-1}\times\CC: r\exp\left({|z'|^2}\right)<|z_n|<R\exp\left({|z'|^2}\right)\Bigr\},
$$
for $0\le r<R\le\infty$. Every $CR$-isomorphism between two $R_{\delta}$-orbits is a composition of a map from $R_{\delta}$ and a dilation in the $z_n$-variable. Therefore, $F$ extends to a holomorphic automorphism of $\CC^n$. Hence, we obtain that $M$ is holomorphically equivalent to $D_r^R$ for some $0\le r<R\le\infty$ and therefore either to $D_{r/R,\,1}$ or (for $R=\infty$) to $D_{0,-1}$.

Suppose now that ${\frak m}=\omega$. In this case ${\cal D}$ is the cylinder ${\cal C}:=\{(z',z_n)\in\CC^{n-1}\times\CC: |z'|<1\}$. The $R_{\omega}$-orbit of every point in ${\cal C}$ is of the form
$$
\left\{(z',z_n)\in\CC^{n-1}\times\CC: |z'|^2+r\exp\left(\hbox{Re}\,z_n\right)=1\right\},
$$
where $r>0$, and any $R_{\omega}$-invariant domain in ${\cal C}$ is of the form
$$
\begin{array}{lll}
\Omega_r^R:&=&\Bigl\{(z',z_n)\in\CC^{n-1}\times\CC: |z'|<1,\\
&&r(1-|z'|^2)<\exp\left(\hbox{Re}\,z_n\right)<R(1-|z'|^2)\Bigr\},
\end{array}
$$
for $0\le r<R\le\infty$. Since every $CR$-isomorphism between two $R_{\omega}$-orbits in ${\cal C}$ is a composition of a map from $R_{\omega}$ and a translation in the $z_n$-variable, $F$ extends to a holomorphic automorphism of ${\cal C}$. In this case $M$ is holomorphically equivalent to $\Omega_r^R$ for some $0\le r<R\le\infty$, and hence either to $\Omega_{r/R,\,1}$ or (for $R=\infty$) to $\Omega_{0,-1}$.

Assume now that ${\frak m}=\varepsilon_{\alpha}$ for some $\alpha>0$. Here ${\cal D}$ is the domain ${\cal C}':={\cal C}\setminus\{z_n=0\}$. 
The $R_{\varepsilon_{\alpha}}$-orbit of every point in ${\cal C}'$ is of the form
$$
\left\{(z',z_n)\in\CC^{n-1}\times\CC: |z'|^2+r|z_n|^{\alpha}=1,\, z_n\ne 0\right\},
$$
where $r>0$, and every $R_{\varepsilon_{\alpha}}$-invariant domain in ${\cal C}'$ is given by
$$
\begin{array}{lll}
{\cal E}_{r,\alpha}^R:&=&\Bigl\{(z',z_n)\in\CC^{n-1}\times\CC: |z'|<1,\\
&&r(1-|z'|^2)^{1/\alpha}<|z_n|<R(1-|z'|^2)^{1/\alpha}\Bigr\},
\end{array}
$$
for $0\le r<R\le\infty$. Since every $CR$-isomorphism between $R_{\varepsilon_{\alpha}}$-orbits is a composition of an element of $R_{\varepsilon_{\alpha}}$ and a dilation in the $z_n$-variable, the map $F$ extends to an automorphism of ${\cal C}'$. Thus, we have shown that $M$ is holomorphically equivalent to ${\cal E}_{r,\alpha}^R$ for some $0\le r<R\le\infty$, and hence either to ${\cal E}_{r/R,1/\alpha}$ or (for $R=\infty$) to ${\cal E}_{0,-1/\alpha}$. 

Suppose that ${\frak m}=B^{n-1}\times\RR$. Here ${\cal D}={\cal C}$, the $R_{B^{n-1}\times\RR}$-orbit of every point in ${\cal C}$ is of the form
$$
b_r:=\left\{(z',z_n)\in\CC^{n-1}\times\CC: |z'|<1,\,\hbox{Im}\,z_n=r\right\}, 
$$
for $r\in\RR$, and any $R_{B^{n-1}\times\RR}$-invariant domain in ${\cal C}$ is given by
$$
B_r^R:=\left\{(z',z_n)\in\CC^{n-1}\times\CC: |z'|<1,\,r<\hbox{Im}\,z_n< R\right\},
$$
for $-\infty\le r<R\le\infty$. We will now show that $F$ extends to an automorphism of ${\cal C}$. Since $F$ maps $b_{r_1}$ onto $b_{r_2}$, for some $r_1,r_2\in\RR$, it has the form $F=\nu\circ g$, where $\nu$ is an imaginary translation in the $z_n$-variable, and $g\in \hbox{Aut}_{CR}(b_{r_1})$. Each of the maps $f$ and $f_1$ transforms the group $G(M)$ into the group $R_{B^{n-1}\times\RR}$, and therefore $g$ lies in the normalizer of $R_{B^{n-1}\times\RR}$ in $\hbox{Aut}_{CR}(b_{r_1})$. Considering $g$ in the general form (\ref{prodgeneral}), we obtain $a_{u_1} a a_{u_1}^{-1}=a_{u_2} a a_{u_2}^{-1}$ for all $a\in\hbox{Aut}(B^{n-1})$ and all $u_1, u_2$. Therefore, $a_{u_1}a_{u_2}^{-1}$ lies in the center of $\hbox{Aut}(B^{n-1})$, which is trivial. Hence we obtain that $a_{u_1}=a_{u_2}$ for all $u_1,u_2$. In addition, there exists $k\in\RR^*$ such that $\mu^{-1}(u)+b\equiv\mu^{-1}(u+kb)$ for all $b\in\RR$. Differentiating this identity with respect to $b$ at $b=0$ we see that $\mu^{-1}(u)=u/k+t$ for some $t\in\RR$. Therefore, $F$ extends to a holomorphic automorphism of ${\cal C}$. This shows that $M$ is holomorphically equivalent to $B_r^R$ for some $0\le r<R\le\infty$, and hence to $B^{n-1}\times\Delta$, which is impossible. Thus, $M$ in fact does not contain orbits with model $B^{n-1}\times\RR$.  

Suppose finally that ${\frak m}=B^{n-1}\times S^1$. Here ${\cal D}={\cal C}'$, the $R_{B^{n-1}\times S^1}$-orbit of every point in ${\cal C}'$ has the form
$$
c_r:=\left\{(z',z_n)\in\CC^{n-1}\times\CC: |z'|<1,\,|z_n|=r\right\},
$$
for $r>0$, and any $R_{B^{n-1}\times S^1}$-invariant domain in ${\cal C}'$ is given by
$$
C_r^R:=\left\{(z',z_n)\in\CC^{n-1}\times\CC: |z'|<1,\,r<|z_n|<R\right\},
$$
for $0\le r<R\le\infty$. We will now show that $F$ extends to an automorphism of ${\cal C}'$. Since $F$ maps $c_{r_1}$ onto $c_{r_2}$, for some $r_1,r_2>0$, it has the form $F=\nu\circ g$, where $\nu$ is a dilation in the $z_n$-variable, and $g\in \hbox{Aut}_{CR}(c_{r_1})$. As in the previous case, each of $f$ and $f_1$ transforms the group $G(M)$ into the group $R_{B^{n-1}\times S^1}$, hence the element $g$ lies in the normalizer of $R_{B^{n-1}\times S^1}$ in $\hbox{Aut}_{CR}(c_{r_1})$. As before, we take $g$ in the general form (\ref{prodgeneral}) and obtain that $a_{u_1}=a_{u_2}$ for all $u_1,u_2$. Furthermore, there exists $k\in\RR^*$ such that $e^{ic}\mu^{-1}(u)\equiv\mu^{-1}(e^{ikc}u)$. Differentiating this identity with respect to $c$ at $c=0$ we see that $\mu^{-1}(u)=e^{it}u^{1/k}$, for some $t\in\RR$. Therefore, $F$ extends to a holomorphic automorphism of ${\cal C}'$. This shows that $M$ is holomorphically equivalent to $C_r^R$ for some $0\le r<R\le\infty$, and hence either to ${\cal E}_{r/R,0}$ or (for $R=\infty$) to ${\cal E}_{0,0}$. This completes the case when $M$ does not contain complex hypersurface orbits.

We will now assume that a complex hypersurface orbit is present in $M$. Recall that there are at most two such orbits and that $G(M)$ can be factored by a normal subgroup isomorphic to $U_1$ to obtain a group isomorphic to $\hbox{Aut}(B^{n-1})$ (see (iii) of Proposition \ref{main1}). It now follows from Propositions \ref{groupsdeterm} and \ref{propleviflat} that the model for every real hypersurface orbit is either $\varepsilon_{\alpha}$ for some $\alpha>0$, or $B^{n-1}\times S^1$. Let $M'$ be the manifold obtained from $M$ by removing all complex hypersurface orbits. It then follows from the above considerations that $M'$ is holomorphically equivalent to either ${\cal E}_{r,\alpha}^R$ or to $C_r^R$ for some $0\le r<R\le\infty$.

Suppose first that $M'$ is equivalent to ${\cal E}_{r,\alpha}^R$ and let $f:M'\ra {\cal E}_{r,\alpha}^R$ be a biholomorphic map satisfying 
(\ref{equivar}) for all $g\in G(M)$, $q\in M'$ and some isomorphism $\varphi: G(M)\ra R_{\varepsilon_{\alpha}}$. The group $R_{\varepsilon_{\alpha}}$ in fact acts on all of ${\cal C}$, and the orbit of any point in ${\cal C}$ with $z_n=0$ is the complex hypersurface        
$$
c_0:=\left\{(z',z_n)\in\CC^{n-1}\times\CC: |z'|<1,\,z_n=0\right\}.
$$
For a point $s\in{\cal C}$ denote by $J_s$ the isotropy subgroup of $s$ under the action of $R_{\varepsilon_{\alpha}}$. If $s_0\in c_0$ and $s_0=(z'_0,0)$, $J_{s_0}$ is isomorphic to $U_1\times U_{n-1}$ and consists of all maps of the form (\ref{autgrp}) with $\theta=1/\alpha$ for which the transformations in the $z'$-variables form the isotropy subgroup of the point $z_0'$ in $\hbox{Aut}(B^{n-1})$.     

Fix $s_0=(z_0',0)\in c_0$ and let
$$
N_{s_0}:=\left\{s\in {\cal E}_{r,\alpha}^R:J_s\subset J_{s_0}\right\}.
$$
We have
$$
N_{s_0}=\left\{(z',z_n)\in\CC^{n-1}\times\CC: z'=z'_0,\,r(1-|z_0'|^2)^{1/\alpha}<|z_n|<R(1-|z_0'|^2)^{1/\alpha}\right\}.
$$
Thus, $N_{s_0}$ is either an annulus (possibly, with an infinite outer radius) or a punctured disk. In particular, $N_{s_0}$ is a complex curve in ${\cal C}'$.

Since $J_{s_0}$ is a maximal compact subgroup of $R_{\varepsilon_{\alpha}}$, $\varphi^{-1}(J_{s_0})$ is a maximal compact subgroup of $G(M)$. Let $O$ be a complex hypersurface orbit in $M$. For $q\in O$ the isotropy subgroup $I_q$ is isomorphic to $U_1\times U_{n-1}$ and therefore is a maximal compact subgroup of $G(M)$ as well. Thus, $\varphi^{-1}(J_{s_0})$ is conjugate to $I_q$ for every $q\in O$ and hence there exists $q_0\in O$ such that $\varphi^{-1}(J_{s_0})=I_{q_0}$. Since the isotropy subgroups in $R_{\varepsilon_{\alpha}}$ of distinct points in $c_0$ do not coincide, such a point $q_0$ is unique.

Let
$$
K_{q_0}:=\left\{q\in M': I_q\subset I_{q_0}\right\}.  
$$
Clearly, $K_{q_0}=f^{-1}(N_{s_0})$. Thus, $K_{q_0}$ is a $I_{q_0}$-invariant complex curve in $M'$ equivalent to either an annulus or a punctured disk. By Bochner's theorem there exist a local holomorphic 
change of coordinates $F$ near $q_0$ on $M$ that identifies an $I_{q_0}$-invariant neighborhood $U$ of $q_0$ with an $L_{q_0}$-invariant neighborhood of the origin in $T_{q_0}(M)$ such that $F(q_0)=0$ and $F(gq)=\alpha_{q_0}(g)F(q)$ for all $g\in I_{q_0}$ and $q\in U$ (here $L_{q_0}$ is the linear isotropy group and $\alpha_{q_0}$ is the isotropy representation at $q_0$). In the proof of Proposition \ref{main1} (see Case 1) we have seen that $L_{q_0}$ has two invariant subspaces in $T_{q_0}(M)$. One of them corresponds in our coordinates to $O$, the other to a complex curve $C$ intersecting $O$ at $q_0$. Observe that near $q_0$ the curve $C$ coincides with $K_{q_0}\cup\{q_0\}$. Therefore, in a neighborhood of $q_0$ the curve $K_{q_0}$ is a punctured analytic disk. Further, if a sequence $\{q_n\}$ from $K_{q_0}$ accumulates to $q_0$, the sequence $\{f(q_n)\}$ accumulates to one of the two ends of $N_{s_0}$, and therefore we have either $r=0$ or $R=\infty$. Since both these conditions cannot be satisfied simultaneously due to hyperbolicity of $M$, we conclude that $O$ is the only complex hypersurface orbit in $M$.

Assume first that $r=0$. We will extend $f$ to a map from $M$ onto the domain
\begin{equation}
\left\{(z',z_n)\in\CC^{n-1}\times\CC: |z'|^2+\frac{1}{R}|z_n|^{\alpha}<1\right\}\label{ellr}
\end{equation}
by setting $f(q_0)=s_0$, where $q_0\in O$ and $s_0\in c_0$ are related as specified above. The extended map is one-to-one and satisfies (\ref{equivar}) for all $g\in G(M)$, $q\in M$. To prove that $f$ is holomorphic on all of $M$, it suffices to show that $f$ is continuous on $O$. It will be more convenient for us to show that $f^{-1}$ is continuous on $c_0$. Let first $\{s_j\}$ be a sequence of points in $c_0$ converging to $s_0$. Then there exists a sequence $\{g_j\}$ of elements of $R_{\varepsilon_{\alpha}}$ converging to the identity such that $s_j=g_js_0$ for all $j$. Then $f^{-1}(s_j)=\varphi^{-1}(g_j)q_0$, and, since $\left\{\varphi^{-1}(g_j)\right\}$ converges to the identity, we obtain that $\{f^{-1}(s_j)\}$ converges to $q_0$. Next, let $\{s_j\}$ be a sequence of points in ${\cal E}_{0,\alpha}^R$ converging to $s_0$. Then we can find a sequence $\{g_j\}$ of elements of $R_{\varepsilon_{\alpha}}$ converging to the identity such that $g_js_j\in N_{s_0}$ for all $j$. Clearly, the sequence $\{f^{-1}(g_js_j)\}$ converges to $q_0$, and hence the sequence $\{f^{-1}(s_j)\}$ converges to $q_0$ as well. Thus, we have shown that $M$ is holomorphically equivalent to domain (\ref{ellr}) and hence to the domain $E_{\alpha}$.

Assume now that $R=\infty$. Observe that the action of the group $R_{\varepsilon_{\alpha}}$ on ${\cal C}$ extends to an action on $\tilde{\cal C}:=B^{n-1}\times\CC\PP^1$ by holomorphic transformations by setting $g(z',\infty):=(a(z'),\infty)$ for every $g\in R_{\varepsilon_{\alpha}}$, where $a$ is the corresponding automorphism of $B^{n-1}$ in the $z'$-variables (see formula (\ref{autgrp})). Now arguing as in the case $r=0$, we can extend $f$ to a biholomorphic map between $M$ and the domain in $\tilde{\cal C}$
$$
\left\{(z',z_n)\in\CC^{n-1}\times\CC: |z'|<1,\,|z_n|>r(1-|z'|^2)^{1/\alpha}\right\}\cup \left(B^{n-1}\times\{\infty\}\right).
$$
This domain is holomorphically equivalent to ${\cal E}_{-1/\alpha}$, and so is $M$.

In the case when $M'$ is holomorphically equivalent to $C_r^R$ for some $0\le r<R\le\infty$, the same argument gives that $M$ has to be equivalent to $B^{n-1}\times\Delta$, which is impossible.  

It now remains to show that all manifolds in (i)-(vii) of (\ref{list}) are pairwise holomorphically non-equivalent. Since the automorphism groups of most manifolds are non-isomorphic (see the discussion following the formulation of Theorem \ref{main} in the introduction), and the orbits of $U_n/\ZZ_m$ in $\CC^n\setminus\{0\}/\ZZ_m$ are topologically different for distinct $m$, we must only prove pairwise non-equivalence of domains within each of the following families: $\left\{\hbox{(ii)-(iv)}\right\}$, $\left\{\hbox{(v)}\right\}$, $\left\{\hbox{(vi)}\right\}$. The first two families consist of Reinhardt domains. It is shown in \cite{Kru} that two hyperbolic Reinhardt domains are holomorphically equivalent if and only if they are equivalent by means of an algebraic map, that is, a map of the form
$$
z_i\mapsto\lambda_i z_1^{a_{i\,1}}\cdot\dots\cdot z_n^{a_{i\,n}},\quad i=1,\dots,n,
$$
where $\lambda_i\in\CC^*$, $a_{ij}\in\ZZ$ for all $i,j$, and $\hbox{det}\,(a_{ij})=\pm 1$. It is straightforward to verify, however, that no two domains within the first two families are algebraically equivalent.

Next, for every domain in the third family its group of holomorphic automorphisms is the group $R_{\omega}$. Therefore, a biholomorphic map between any two domains $D_1$ and $D_2$ in this family, takes every $R_{\omega}$-orbit from $D_1$ into a $R_{\omega}$-orbit in $D_2$. However, as we noted above, every $CR$-isomorphism between two $R_{\omega}$-orbits is a composition of an element of $R_{\omega}$ and a translation in the $z_n$-variable. Therefore, $D_1$ can be mapped onto $D_2$ by such a translation. It is clear, however, that no two domains in the third family can be obtained from one another by translating the $z_n$-variable.

This completes the proof of Theorem \ref{main} in the non-homogeneous case.

\section{The Homogeneous Case}\label{proofmain2}

\setcounter{equation}{0}

If $M$ is homogeneous, by \cite{N}, \cite{P-S} it is holomorphically equivalent to a Siegel domain of the second kind in $\CC^n$. For $n=2$, this gives that $M$ is equivalent to either $B^2$ or $\Delta^2$, which is impossible since $d(B^2)=8$ and $d(\Delta^2)=6$. For $n=3$ we obtain that $M$ is equivalent to one of the following domains: $B^3$, $B^2\times\Delta$, $\Delta^3$, $S$, where $S$ is the 3-dimensional Siegel space mentioned in the introduction. Among these domains only $\Delta^3$ has an automorphism group of dimension 9.

For $n\ge 4$ we will prove the following proposition.

\begin{proposition}\label{siegel}\sl Let $U\subset\CC^n$, $n\ge 4$, be a Siegel domain of the second kind. Suppose that $d(U)=n^2$. Then $n=4$ and $U$ is holomorphically equivalent to $B^2\times B^2$.
\end{proposition}

\noindent {\bf Proof:} Our proof is similar to that of Proposition 3.1 of \cite{IKra}. The domain $U$ has the form
$$
U=\left\{(z,w)\in\CC^{n-k}\times\CC^k: \hbox{Im}\,w-F(z,z)\in C\right\},
$$
where $1\le k\le n$, $C$ is an open convex cone in $\RR^k$ not containing an entire affine line and $F=(F_1,\dots,F_k)$ is a $\CC^k$-valued Hermitian form on $\CC^{n-k}\times\CC^{n-k}$ such that $F(z,z)\in\overline{C}\setminus\{0\}$ for all non-zero $z\in\CC^{n-k}$.

We will first show that $k\le 2$. As we noted in \cite{IKra}
\begin{equation}
d(U)\le 4n-2k+\hbox{dim}\,{\frak g}_0(U).\label{form1}
\end{equation}
Here ${\frak g}_0(U)$ is the Lie algebra of all vector fields on $\CC^n$ of the form
$$
X_{A,B}=Az\frac{\partial}{\partial z}+Bw\frac{\partial}{\partial w},
$$
where $A\in{\frak{gl}}_{n-k}(\CC)$, $B$ belongs to the Lie algebra ${\frak g}(C)$ of the group of linear automorphisms of the cone $C$, and the following holds:
\begin{equation}
F(Az,z)+F(z,Az)=BF(z,z),\label{form3}
\end{equation}
for all $z\in\CC^{n-k}$.

In \cite{IKra} we showed that
$$
\hbox{dim}\,{\frak g}_0(U)\le\frac{3k^2}{2}-k\left(2n+\frac{1}{2}\right)+n^2+1,
$$ 
which together with (\ref{form1}) gives
\begin{equation}
d(U)\le\frac{3k^2}{2}-k\left(2n+\frac{5}{2}\right)+n^2+4n+1.\label{form2}
\end{equation}
It is straightforward to check that the right-hand side of (\ref{form2}) is strictly less than $n^2$ if $k\ge 3$, hence $k\le 2$.

If $k=1$, the domain $U$ is equivalent to $B^n$ which is impossible. Hence $k=2$. It follows from (\ref{form3}) that the matrix $A$ is determined by the matrix $B$ up to a matrix $L\in{\frak {gl}}_{n-2}(\CC)$ satisfying
$$
F(Lz,z)+F(z,Lz)=0,
$$
for all $z\in\CC^{n-2}$. Let $s$ be the dimension of the subspace of all such matrices $L$. Then
$$
\hbox{dim}\,{\frak g}_0(U)\le s+\hbox{dim}\,{\frak g}(C),
$$
and Lemma 3.2 of \cite{IKra} yields
$$
\hbox{dim}\,{\frak g}_0(U)\le s+2,
$$
which, together with (\ref{form1}) implies
\begin{equation}
s\ge n^2-4n+2.\label{form4}
\end{equation}
 
By the definition of Siegel domain, there exists a positive-definite linear combination of the components of $F$, and we can assume that $F_1$ is positive-definite. Further, applying an appropriate linear transformation of the $z$-variables, we can assume that $F_1$ is given by the identity matrix and $F_2$ by a diagonal matrix.

Suppose first that the matrix of $F_2$ is scalar. If $F_2\equiv 0$, then $U$ is holomorphically equivalent to $B^{n-1}\times\Delta$ which is impossible. If $F_2\not\equiv 0$, then $U$ is holomorphically equivalent to the domain
$$
V:=\left\{(z,w)\in\CC^{n-2}\times\CC^2:\hbox{Im}\, w_1-|z|^2>0,\,\hbox{Im}\,w_2-|z|^2>0\right\}.
$$
It was shown in \cite{IKra} that $d(V)\le n^2-2n+3$ and hence $d(V)<n^2$. Thus, the matrix of $F_2$ is not scalar. Inequality (\ref{form4}) now yields that the matrix of $F_2$ can have at most one pair of distinct eigenvalues. Therefore $n=4$ and thus $U$ is holomorphically equivalent to $B^2\times B^2$.

The proof of Proposition \ref{siegel} is complete.\qed
\smallskip\\

Theorem \ref{main} in the homogeneous case now follows from Proposition \ref{siegel} and the preceding remarks.\qed

\section{Examples for the case $d(M)=n^2-2$}\label{examples}

\setcounter{equation}{0}

In this section we will give an example of a family of pairwise holomorphically non-equivalent smoothly bounded Reinhardt domains in $\CC^2$ with automorphism group of dimension 2. This family is parametrized by sets in $\RR^2$ that satisfy only very mild conditions. In particular, no explicit formulas describe the domains in the family.  

Choose a set $Q\subset\RR^2_{+}:=\{(x,y)\in\RR^2: x\ge 0, y\ge 0\}$ in such a way that the associated set in $\CC^2$
$$
D_Q:=\left\{(z_1,z_2)\in\CC^2: (|z_1|, |z_2|)\in Q\right\}
$$
is a smoothly bounded Reinhardt domain and contains the origin. By
\cite{Su}, two bounded Reinhardt domains containing the origin are
holomorphically equivalent if and only if one is obtained from the other
by means of a dilation and permutation of coordinates. In \cite{FIK}
all smoothly bounded Reinhardt domains with non-compact automorphism group were listed. All these domains contain the origin, and it is not difficult to choose $Q$ so that $D_Q$ is not holomorphically equivalent to any of the domains from \cite{FIK} and thus ensure that $\hbox{Aut}(D_Q)$ is compact. It then follows from the explicit description of the automorphism groups of bounded Reinhardt domains (see \cite{Kru}, \cite{Sh}) that $\hbox{Aut}(D_Q)^c$ is isomorphic to either $U_2$ or $U_1\times U_1$. Theorem 1.9 of \cite{GIK} (or, alternatively, Theorem \ref{main} above) now yields that there does not exist a hyperbolic Reinhardt domain in $\CC^2$ containing the origin for which the automorphism group is four-dimensional. Hence $\hbox{Aut}(D_Q)^c$ is in fact isomorphic to $U_1\times U_1$ and thus $d(D_Q)=2$. 

The freedom in choosing a set $Q$ that satisfies the above requirements is very substantial, and by varying $Q$ one can produce a family of pairwise non-equivalent Reinhardt domains $D_Q$ that cannot be described by explicit formulas.

{\obeylines
Department of Mathematics
The Australian National University
Canberra, ACT 0200
AUSTRALIA
E-mail: alexander.isaev@maths.anu.edu.au
}

\end{document}